\documentclass[a4paper]{amsart}

\usepackage{amssymb}
\usepackage{amsmath,amsthm}
\usepackage{enumerate,paralist}
\usepackage{amstext}
\usepackage{psfrag}
\usepackage{dsfont}
\usepackage[dvips]{epsfig}
\usepackage[english]{babel}
\usepackage{cite}
\usepackage{hyperref}
\usepackage{color}
\usepackage{cleveref}
\usepackage[svgnames]{xcolor}
\hypersetup{hidelinks,colorlinks=true,allcolors=DarkBlue}

\theoremstyle{plain}
\newtheorem{theorem}{Theorem}

\theoremstyle{definition}

\DeclareMathOperator{\sign}{sign}

\def\R{{\mathbb R}}

\def\Z{{\mathbb Z}}

\def\text#1{\mbox{\rm #1\,}}

\begin{document}
\frenchspacing

\title[Asymptotically optimal control]
{Asymptotically optimal control \\ for a simplest distributed system}

\author{Alexander Ovseevich}
\address
{
Institute for Problems in Mechanics, Russian Academy of Sciences \\ 119526, Vernadsky av., 101/1, Moscow, Russia \\
ovseev@ipmnet.ru 
}

\author{Aleksey Fedorov}
\address
{
LPTMS, CNRS, Univ. Paris-Sud, Universit\'e Paris-Saclay, 91405, 100 B\^atiment, 15 rue Georges Cl\'emenceau, Orsay, France \\
aleksey.fedorov@lptms.u-psud.fr
}

\maketitle

\begin{abstract}
We study the problem of the minimum-time damping of a closed string under a bounded load, applied at a single fixed point. 
A constructive feedback control law is designed, which allows bringing the system to a bounded neighbourhood of the terminal manifold. 
The law has the form of the dry friction at the point, where the load is applied. 
The motion under the control is governed by a nonlinear wave equation. 
The existence and uniqueness of solution of the Cauchy problem for this equation are proved. 
The main result is the asymptotic optimality of the suggested control law.

\medskip\noindent
\textsc{Keywords} maximum principle, reachable sets, linear system

\medskip\noindent
\textsc{MSC 2010:} 93B03, 93B07, 93B52.
\end{abstract}

\section{Introduction}

Oscillating systems are one of the central objects of control theory~\cite{Akulenko,Lions}. 
In particular, problems of control for oscillating distributed systems are important for a wide spectrum of technological objects~\cite{Butkovskiy}.
At present in practice there are essentially no design methods of constructive control laws, which combine the implementation simplicity with the optimality properties. 
One of the approaches to design of such control laws is the use of the asymptotic theory of reachable sets~\cite{Ovseevich}. 
Based on this theory, a method of control for a system of an arbitrary number of linear oscillators was proposed~\cite{Ovseevich2}. 
In the present work, the technique developed in~\cite{Ovseevich2,Ovseevich3} is applied to the design of an asymptotically optimal control for a string. 
The main result is the construction of a bounded feedback control possessing an asymptotic optimality: 
the ratio of the time of steering to a fixed bounded neighborhood of the terminal manifold under this control to the minimum-time is close to 1,
if the initial energy is large.

\section{Problem statement}

The controlled object is the closed string described by the following equation:
\begin{equation}\label{string_eq}
	\frac{\partial^2 f}{\partial t^2}=\frac{\partial^2 f}{\partial x^2}+u\delta, \quad |u|\leq1,
\end{equation}
where the spatial variable $x$ (angle) runs over the one-dimensional torus $\mathcal{T}=\mathbb{R}/2\pi\mathbb{Z}$, 
$t$ is time, 
and $\delta=\delta(x)$ denotes the Dirac $\delta$-function. 
The phase space $S$ of the control system consists of pairs $\mathfrak{f}=(f_0,f_1)$ of even distributions on $\mathcal{T}$. 
Thus, the systems has the form
\begin{equation}\label{string_eq2}
	\frac{\partial \mathfrak{f}}{\partial t}=A\mathfrak{f}+Bu,
	\quad
	|u|\leq1,
\end{equation}
where
\begin{equation}
	A=\left(
	\begin{array}{cc}
	0 & 1 \\
	\Delta & 0 \\
	\end{array}
	\right),
	\quad
	\Delta=\frac{\partial^2 }{\partial x^2},
	\quad
	B=\left(
	\begin{array}{c}
	0 \\
	\delta\\
	\end{array}
	\right).
\end{equation}

Parity can be explained by the fact that any solution of Eq. (\ref{string_eq}) with zero initial conditions is even.
Thus, the odd part of possible states of the string cannot be steered to the terminal set.

A mechanical model of such system can be not only a string, but also, e.g., a toroidal acoustic resonator.

The string equation (\ref{string_eq}) can be rewritten in the form
\begin{equation}\label{model}
	\left(\frac{\partial }{\partial t}-\frac{\partial}{\partial x}\right)g(x,t)=\delta(x)u(t),
	\quad
	|u|\leq1,
\end{equation}
where
\begin{equation}\label{g}
	g=\frac{\partial f_0}{\partial x}+f_1.
\end{equation}
The knowledge of the function $g$ is equivalent to the knowledge of ${\partial f_0}/{\partial x}$ and $f_1$, 
because the function ${\partial f_0}/{\partial x}$ is odd, $f_1$ is even. 
The derivative ${\partial f_0}/{\partial x}$ defines $f_0$ up to an additive constant, 
which is irrelevant if the goal of the control is to stop the oscillations or to stop the motion at an unspecified point.

Our goal is to design an easily implementable feedback control providing the fastest transfer from the initial state to the terminal manifold ${\mathcal C}$. 
The terminal manifold consists of pairs of constants $\{(c_0,c_1)^*\in\R^2\subset S\}$. 
In the present work, we consider three problems:
\begin{enumerate}
	\item Complete stop at a given point: ${\mathcal C}=0$.
	\item Stop moving: ${\mathcal C}=\R\times 0$.
	\item Oscillation damping: ${\mathcal C}=\R^2$.
\end{enumerate}
We can regard the factor space $\overline{S}=S/{\mathcal C}$ as the phase space, since the terminal manifold ${\mathcal C}$  is invariant wrt the phase flow of system (\ref{string_eq2}). 
In the present paper, we put an accent on Problem 2 and Problem 3, where the string model (\ref{model})--(\ref{g}) is applicable.

\section{Reachable sets}

The controllability is the first important issue in construction of the control for the system. 
We study controllability by computing reachable sets. 
This approach has much in common with the methods suggested in~\cite{Lions}.

By passing to the backward time, we reduce the problem of characterization of states allowing damping to computation of the reachable set of system (\ref{string_eq2}) from zero. 
The reachable set $D(T)$, $T\geq 0$ of system (\ref{string_eq2}) is a closed convex set. 
Computation of the support function $H=H_{D(T)}(\xi)$, where $\xi\in\overline{S}^*$ is the vector of distributions orthogonal to ${\mathcal C}$, gives:
\begin{equation}\label{support2}
	{D(T)}(\xi)=\int_0^T\left|\sum_{n\neq0}\left(\psi_n\cos nt+\frac{\phi_n}{n}\sin nt\right)+\psi_0+\phi_0 t\right|dt=\int_0^{T}\left|\zeta(t)\right|dt.
\end{equation}
Here
\begin{equation}\label{cosine}
	\zeta(t)=\xi_1(t)+\int_0^t\xi_0(x)dx, \mbox { where } \xi_0(x)=\sum{\phi_n\cos{nx}}, \,\,\, \xi_1(x)=\sum{\psi_n\cos{nx}}
\end{equation}
is the Fourier expansion of the even functions $\xi_0$, and $\xi_1$.

Let $\mathcal{D}$ be the set of vectors reachable in an arbitrary time. 
Then, the following characterization of states allowing damping holds true:
\begin{theorem}\label{D}
	The set $\mathcal{D}$ consists precisely of pairs $(f_0,f_1)$, where $f_0$ is (an even) Lipschitz function and $f_1$ is (an even) bounded function. 
	This means that function (\ref{g}) $g=\frac{\partial f_0}{\partial x}+f_1$ is bounded.
\end{theorem}
In other words, the reachability criteria is that
\begin{equation}\label{reach2}
	\rho(\mathfrak{f}):=\left\Vert\frac{\partial f_0}{\partial x}+f_1\right\Vert=
	\sup_{x\in\mathcal{T}}\left|\left(\frac{\partial f_0}{\partial x}+f_1\right)(x)\right|<\infty.
\end{equation}

Consider problems of the complete stop at an unspecified point and of the oscillation damping. 
The corresponding reachable sets $\overline{D}(T)$ belong to the reduced phase spaces $\overline{S}=S/{\mathcal C}$. 
In this case, we have the following asymptotic result:
\begin{theorem}\label{32}
As $T\to\infty$ the following limit formula holds true:
\begin{equation}\label{asympt3}
	\frac1TH_{\overline{D}(T)}(\xi)\rightarrow\hfill\frac1{2\pi}\int_0^{2\pi}\left|\sum\left(\psi_n\cos nt+\frac{\phi_n}{n}\sin nt\right)\right|dt=\frac1{2\pi}\int_0^{2\pi}\left|\zeta(t)\right|dt.
\end{equation}
where $\zeta(t)=\xi_1(t)+\int_0^t\xi_0(x)dx$.
\end{theorem}
One can put the equal sign instead of the limit sign provided that $T=2k\pi$, $k\in\Z$. 
The case $\psi_0=0$ corresponds to the oscillation damping problem. 
The equality $\phi_0=0$ holds in both problems 2 and 3.

Define a convex closed subset $\Omega$ in the phase space in such a way that the support function $H_\Omega(\xi)$ is given by formula (\ref{asympt3}):
\begin{theorem}\label{omega}
	In the problem of the complete stop at an unspecified point the convex set $\Omega$ is defined as follows:
	$$
		\Omega=\left\{\mathfrak{f}=(f_0,f_1): \rho(\mathfrak{f})={2\pi}\sup_{x\in\mathcal{T}}\left\vert\frac{\partial f_0}{\partial x}(x)+f_1(x)\right\vert\leq1\right\}.
	$$
	In the oscillation damping problem
	$$
		\Omega=\left\{\mathfrak{f}=(f_0,f_1): \rho(\mathfrak{f})={2\pi}\inf_{c\in\R}\sup_{x\in\mathcal{T}}\left\vert\frac{\partial f_0}{\partial x}(x)+f_1(x)+c\right\vert\leq1\right\}.
	$$
\end{theorem}
Consider the problem of the complete stop at an unspecified point or the oscillation damping problem
\begin{theorem}\label{rho00}\label{time}
	Let $T$ be the optimal time of the transfer from the initial state $\mathfrak{f}$ to the target manifold ${\mathcal C}$.
	Then ${T}\thicksim{\rho(\mathfrak{f})}$ as $\rho(\mathfrak{f})\to\infty$.
\end{theorem}
This statement follows from Theorem \ref{32} and Theorem \ref{omega}.

\section{Control in the generalized dry-friction form}

In the stop moving problem we design a ``quasioptimal'' control $u(\mathfrak{f})$ based on the following idea: 
The optimal control in the state ${\mathfrak f}$ implements the steepest decent in the direction orthogonal to the boundary of the reachable set $D(T)$ passing through $\mathfrak{f}$. 
The control in the generalized dry-friction form in problems of stop moving and oscillations damping 
implements the steepest decent in the direction orthogonal to the boundary of the approximate reachable set $T\Omega$. 
This means that
\begin{equation}\label{Legendre2}
	u({\mathfrak f})=-\sign\langle B,\xi\rangle=-\sign\xi_1(0),=-\sign\zeta(0),
\end{equation}
where the momentum $\xi$ is to be found from the equation
\begin{eqnarray}\label{asympt5}
	T^{-1}{\mathfrak f}=\frac{\partial H_{\Omega}}{\partial\xi}(\xi),
\end{eqnarray}
and $\zeta(t)=\xi_1(t)+\int_0^t\xi_0(x)dx$. 
Such a strategy of the control design was used for the system of an arbitrary number of linear oscillators~\cite{Ovseevich2,Ovseevich3}. 
The design of the control for the string is obtained by the formal passage to the limit $N\to\infty$, 
where $N$ in the number of oscillators in the construction considered in the cited papers.

In the stop moving problem, we have ${\mathfrak f}=(f_0,f_1)$, where
\begin{eqnarray}\label{asympt6}
	T^{-1}f_0(x)=-\int_0^{x}(\sign\zeta(y))^-dy,\quad T^{-1}f_1(x)=(\sign\zeta(x))^+.
\end{eqnarray}
Here $f^\pm$ denotes the even/odd part of the function $f$:
\begin{equation}\label{pm}
	f^\pm(x)=\frac12(f(x)\pm f(-x)).
\end{equation}
Since $(\sign\zeta(0))^+{=}\sign\zeta(0)$, 
we then obtain that $\sign\zeta(0)=\sign f_1(0)$ and on the formal level the motion of the string is described by the nonlinear wave equation:
\begin{equation}\label{string_eq22}
	\frac{\partial^2 f}{\partial t^2}=\frac{\partial^2 f}{\partial x^2}-\sign\left(\frac{\partial f}{\partial t}(0)\right)\delta.
\end{equation}

In the first-order model (\ref{model}), the equations of motion are
\begin{equation}\label{model22}
	\left(\frac{\partial }{\partial t}-\frac{\partial}{\partial x}\right)g(x,t)=-\delta(x)\sign(g(0,t)).
\end{equation}
This is equivalent to
\begin{equation}\label{newmodel5}
	g(z,t)=g(z+t,0)-\sum_{J}\sign g(0,z+t+2k\pi),
\end{equation}
where the set of indices $J={J_t}$ consists of $k\in\Z$ such that $z+2k\pi\in[-t,0]$. 
Eq. (\ref{newmodel5}) should not be understood pointwise, but as the equality of distributions wrt $z$ depending on the parameter $t$.

\section{Existence of the motion under the dry-friction control}

As in the case of a system of an arbitrary number of linear oscillators, the issue of the existence of the motion under the control plays a special role. 
In ~\cite{Ovseevich3}, this question is resolved within the framework of the DiPerna--Lions theory of the Cauchy problem for ordinary differential equations with singularities.

However, the question of existence and uniqueness of the motion for the sting requires different methods.
\begin{theorem}\label{motion}
	The Cauchy problem for the nonlinear wave equation
	\begin{equation}\label{model_theorem}
		\left(\frac{\partial }{\partial t}-\frac{\partial}{\partial x}\right)g(x,t)=-\delta(x)u(t),\quad u(t)=\sign{g(0,t)},
	\end{equation}
	where $g(x,0)$ is a given bounded measurable function, has a unique bounded solution for $t\geq0$. 
	Functions $\phi(t)=g(0,t)$ and $G(x)=g(x,0)$ form the only solution of the equation
	\begin{equation}\label{model501}
		\phi(t)+ \frac12 v(t)+\sum_{k\neq0,2k\pi\in[-t,0]}v(t+2k\pi)=G(t),\quad v=\sign\phi,
	\end{equation}
	where $G$ is the given function, $v$ and $\phi$ are unknown. 
	In Eq.~(\ref{model501}), the function $\sign(x)$ is interpreted as multivalued, meaning that $\sign(0)=[-1,1]$.
\end{theorem}

Note that solution of Eq.~(\ref{model501}), and thus the non-linear wave equation~(\ref{model_theorem}), is reached by algebraic operations {\it only}. 
If $t<2\pi$, Eq.~(\ref{model501}) reduces to
\begin{equation}\label{model502}
	\phi(t)+ \frac12v(t)=G(t),\quad v={\rm sign}\phi,
\end{equation}
which has a unique solution $\phi(t)=G(t)-\frac12$ at $G(t)>\frac12$, and $\phi(t)=G(t)+\frac12$ at $G(t)<-\frac12$.
Otherwise, $\phi(t)=0$ and $v(t)=2G(t)$. 
From Eq.~(\ref{model501}) and periodicity of $G(t+2\pi)=G(t)$, we obtain equation
\begin{equation}\label{phi3}
	\phi(t+2\pi)+ \frac12\sign \phi(t+2\pi)=G(t)-\sign \phi(t),
\end{equation}
which allows to extend the function $\phi(t)$ to $t$ in the interval $(0,4\pi)$. 
Similar arguments give an extension to all positive values of $t$.

\section{Asymptotic optimality}
Our main result is as follows:
\begin{theorem}\label{AsymptoticOptimality}
	Consider the evolution $\rho(t)=\rho(\mathfrak{f}_t)$ of the functional $\rho(\mathfrak{f})$, when $\mathfrak{f}$ changes according to the nonlinear wave equation~(\ref{string_eq22}). 
	Put
	\begin{equation}
		M=\min\{\rho(0),\rho(T)\}.
	\end{equation}
	Then, as $M\to+\infty$, $T\to+\infty$, we have
	\begin{equation}\label{approx_T}
		{(\rho(0)-\rho(T))}/{T}=1+O(1/{T}+1/{M}).
	\end{equation}
	For any other admissible control,
	\begin{equation}\label{approx_T2}
		{(\rho(0)-\rho(T))}/{T}\leq 1+O(1/{T}+1/{M}).
	\end{equation}
\end{theorem}

According to Theorem \ref{time}, the functional $\rho(\mathfrak{f})$ is asymptotically close to the optimal time of passage from $\mathfrak{f}$ to the terminal manifold. 
Therefore, the result of Theorem \ref{AsymptoticOptimality} says that along a trajectory of the motion under the generalised dry-friction control, 
the passage time behaves in the same way as along a time-optimal trajectory.

Thus, for the simplest distributed system the control in the generalized dry-friction form is asymptotically time-optimal in a precise sense.

\bigskip
\noindent{\bf Acknowledgements}. This work was supported by the Russian Science Foundation, grant 16-11-10343.

\end{document}